\documentclass[a4paper,11pt,oneside]{article}

\usepackage[T1]{fontenc}
\usepackage{amsfonts}
\usepackage{amssymb}
\usepackage{amsmath}

\usepackage{graphicx}
\usepackage{subfigure}

\usepackage{latexsym}

\usepackage{verbatim}

\makeatletter
\renewcommand{\@evenfoot}{\hfil - \thepage\ - \hfil}
\renewcommand{\@oddfoot}{\hfil - \thepage\ - \hfil}
\makeatother

\newcounter{theo}

\newenvironment{thm*}
{\noindent\hangafter=1\hangindent=15pt\refstepcounter{theo}\textsc{Theorem --}\begin{sffamily}}
{\end{sffamily}\par}

\newenvironment{thmt*}[1]
{\noindent\hangafter=1\hangindent=15pt\refstepcounter{theo}\textsc{Theorem -- #1}\\\begin{sffamily}}
{\end{sffamily}\par}

\newenvironment{prop*}
{\noindent\hangafter=1\hangindent=15pt\textsc{Proposition. --}\begin{sffamily}}
{\end{sffamily}\par}

\newenvironment{lemme*}
{\noindent\hangafter=1\hangindent=15pt\refstepcounter{theo}\textsc{Lemma. -- }\begin{sffamily}}
{\end{sffamily}\par}

\newenvironment{proof}
{\noindent\textsc{Proof.}\setlength{\parskip}{2pt}\setlength{\parindent}{0pt}}
{\hfill$\square$\par}

\newcommand{\bP}{\mathbb{P}}

\newcommand{\N}{\mathbb{N}}

\newcommand{\eqd}{\begin{eqnarray*}} 
\newcommand{\eqf}{\end{eqnarray*}} 

\usepackage{hyperref}

\newcommand{\limites}[2]{\overset{#1}{\underset{#2}{\longrightarrow}}}
\title{A note on the recurrence of edge reinforced random walks}
\author{Laurent Tournier\footnote{Universit\'e de Lyon; CNRS; Universit\'e Lyon 1, Intitut Camille Jordan, 43 bld du 11 novembre 1918, F-69622 Villeurbanne Cedex, France. E-mail address: {\tt tournier@math.univ-lyon1.fr}}}
\date{}

\begin{document}

\begin{center}
\LARGE \bf
A note on the recurrence of edge reinforced random walks

\vspace{1cm}
\large
Laurent Tournier\footnote{Universit\'e de Lyon; CNRS; Universit\'e Lyon 1, Intitut Camille Jordan, 43 bld du 11 novembre 1918, F-69622 Villeurbanne Cedex, France. E-mail address: {\tt tournier@math.univ-lyon1.fr}}

\end{center}
\vspace{1cm}

\begin{center}
\begin{minipage}{15cm}
\small \textbf{Abstract.} We give a short proof of Theorem 2.1 from~\cite{MR07}, stating that the linearly edge reinforced random walk (ERRW) on a locally finite graph is recurrent if and only if it returns to its starting point almost surely. This result was proved in~\cite{MR07} by means of the much stronger property that the law of the ERRW is a mixture of Markov chains. Our proof only uses this latter property on \emph{finite} graphs, in which case it is a consequence of De Finetti's theorem on exchangeability. 
\end{minipage}
\end{center}

Although the question of the recurrence of linearly edge reinforced random walks (ERRW) on infinite graphs has known important breakthroughs in the recent years (cf.~notably \cite{MR09}), it seems that the only known proof that one almost-sure return implies the recurrence of the walk is based on the difficult fact that ERRWs on infinite graphs are mixtures of Markov chains (cf.~\cite{MR07}). We provide in this note a short and simple proof of that property, with the finite case as the only tool. 

Let $G=(V,E)$ be a locally finite undirected graph, and $\alpha=(\alpha_e)_{e\in E}$ be a family of positive real numbers. The \emph{linearly edge reinforced random walk on $G$ with initial weights $\alpha$ starting at $o\in V$} is the nearest-neighbour random walk $(X_k)_{k\geq 0}$ on $V$ defined as follows: $X_0=o$; then, at each step, the walk crosses a neighbouring edge chosen with a probability proportional to its weight; and the weight of an edge is increased by 1 after it is traversed. 

The only property to be used in this note is the following consequence of De Finetti's theorem for Markov chains (cf.~\cite{DF80}, and~\cite{KR99} for instance): if $G$ is finite, then there exists a probability measure $\mu$ on transition matrices on $G$ such that the law of the ERRW $X$ is $\int P_\omega(\cdot) d\mu(\omega)$ where $P_\omega$ is the law of the Markov chain on $V$ with transition $\omega$ starting at $o$. 

Here is the statement of the (main) part of Theorem 2.1 in~\cite{MR07} (cf.~remark after the proof). 

\begin{thm*}
For the linearly edge-reinforced random walk (ERRW) on any locally finite weighted graph, the following two statements are equivalent: 
\begin{enumerate}
	\item the ERRW returns to its starting point with probability 1;
	\item the ERRW returns to its starting point infinitely often with probability 1. 
\end{enumerate}
\end{thm*}

\begin{proof}
On finite graphs, this result follows from a Borel-Cantelli argument (cf.~\cite{KR99} and the remark after the proof). 
Let us therefore denote by $\bP$ the law of the ERRW on an \emph{infinite} locally finite weighted graph $G$ starting at $o$. Assume that condition {\sf (i)} holds. 

For any $n\in\N$, we introduce the finite graph $G_n$ defined from the ball $B(n+1)$ of center $o$ and radius $n+1$ in $G$ by identifying the points at distance $n+1$ from $o$ to a new point $\delta_n$. The law of the ERRW on $G_n$ (with same weights as in $G$) starting at $o$ is denoted by $\bP_{G_n}$. 

Let us also define the successive return times $\tau^{(1)},\tau^{(2)},\ldots$ of the ERRW at $o$, the exit time $T_n$ from $B(n)$, and the hitting time $\tau_{\delta_n}$ of $\delta_n$ in $G_n$. Note that the laws $\bP$ and $\bP_{G_n}$ may be naturally coupled in such a way that the trajectories coincide up to time $T_n=\tau_{\delta_n}$. 

We have, for all $k\geq 1$, 
\[\bP(\tau^{(k)}<\infty)
	=  \bP(\tau^{(k)}<\infty, \tau^{(k)}< T_n)+\bP(T_n< \tau^{(k)}<\infty),\]
and the second term converges to 0 when $n\to\infty$ since $T_n\geq n\limites{}{n}\infty$. Therefore, 
\begin{eqnarray}
\bP(\tau^{(k)}<\infty) 
	& = &\bP(\tau^{(k)}<\infty, \tau^{(k)}< T_n)+o_n(1)\notag\\
	& = &\bP_{G_n}(\tau^{(k)}<\infty,\tau^{(k)}<\tau_{\delta_n})+o_n(1)\label{eqn:tau2}.
\end{eqnarray}
(NB: the condition $\tau^{(k)}<\infty$ on last line could be dropped since {\sf (ii)} is true for ERRW on finite graphs). In particular, assumption {\sf(i)} gives: 
\begin{equation}\label{eqn:tau1}
\lim_n\bP_{G_n}(\tau^{(1)}<\infty, \tau^{(1)}<\tau_{\delta_n})=\bP(\tau^{(1)}<\infty)=1.
\end{equation}

Since $G_n$ is finite, we may write $\bP_{G_n}$ as a mixture of Markov chains: $\bP_{G_n}(\cdot)=\int P_{G_n,\omega}(\cdot) d\mu_n(\omega)$. Thus we have, according to (\ref{eqn:tau1}), 
\begin{equation*}
\lim_n \int P_{G_n,\omega}(\tau^{(1)}<\infty, \tau^{(1)}<\tau_{\delta_n}) d\mu_n=1
\end{equation*}
and, for all $k\geq 1$, according to (\ref{eqn:tau2}) and Markov property (applied $k-1$ times),
\begin{eqnarray*}
\bP(\tau^{(k)}<\infty)
	& = & \lim_n \int P_{G_n,\omega}(\tau^{(k)}<\infty, \tau^{(k)}<\tau_{\delta_n})\,d\mu_n\\
	& = & \lim_n \int P_{G_n,\omega}(\tau^{(1)}<\infty, \tau^{(1)}<\tau_{\delta_n})^k\, d\mu_n.
\end{eqnarray*}
We may conclude that the last limit equals 1 thanks to the following very simple Lemma: 

\begin{lemme*}
If $(f_n)_n, (\mu_n)_n$ are respectively a sequences of measurable functions and probability measures such that, for all $n$, $0\leq f_n\leq 1$, and $\int f_n d\mu_n\limites{}{n} 1$, then: $$\mbox{for every integer $k\geq 1$, }\int (f_n)^k d\mu_n \limites{}{n} 1.$$
\end{lemme*}
\begin{proof}
Indeed, we have $0\leq f_n^k\leq 1$, hence: 
$$0\leq 1-\int f_n^k d\mu_n = \int (1-f_n^k) d\mu_n = \int (1-f_n)(1+f_n+\cdots+f_n^{k-1})d\mu_n\leq k\int (1-f_n)d\mu_n\limites{}{n} 0.$$
\end{proof}

As a conclusion, $\bP(\tau^{(k)}<\infty)=1$ for all $k\geq 1$, hence $\bP(\forall k,\tau^{(k)}<\infty)=1$, which is {\sf(ii)}.
\end{proof}

\paragraph{Remark} Condition {\sf(ii)} implies that the ERRW visits every edge in the connected component of the starting point infinitely often in both directions, by means of the conditional Borel-Cantelli lemma, cf.~the end of the proof of Theorem 1.1 in \cite{MR09} or Proposition 1 of~\cite{KR99} for a direct proof. 


\end{document}